\documentclass[12pt]{article}
\usepackage[cp1251]{inputenc}

\usepackage[T2A]{fontenc}

\usepackage[english]{babel}
\usepackage{amsfonts}
\usepackage{enumerate}
\usepackage{amssymb}

\begin{document}

\begin{center}
\textbf{Independent linear forms on the group $\Omega_p$}
\end{center}

\bigskip

\begin{center}
\textbf{Margaryta Myronyuk}
\end{center}

\bigskip

 \begin{abstract}

Let $\Omega_p$ be the group of $p$-adic numbers,  $ \xi_1$, $\xi_2$,
$\xi_3$ be independent random variables with values in $\Omega_p$
and   distributions   $\mu_1$, $\mu_2$, $\mu_3$. Let $\alpha_j,
\beta_j, \gamma_j$ be topological automorphisms of $\Omega_p$. We
consider linear forms $L_1 = \alpha_1\xi_1 + \alpha_2 \xi_2+
\alpha_3 \xi_3$, $L_2=\beta_1\xi_1 + \beta_2 \xi_2+ \beta_3 \xi_3$
and $L_3=\gamma_1\xi_1 + \gamma_2 \xi_2+ \gamma_3 \xi_3$. Assuming
that the linear forms $L_1$, $L_2$ and $L_3$ are independent, we
describe possible distributions $\mu_1$, $\mu_2$, $\mu_3$. This
theorem is an analogue of the well-known Skitovich-Darmois theorem,
where a Gaussian distribution on the real line is characterized by
the independence of two linear forms.
\end{abstract}

\bigskip

\textbf{Keywords}: Group of p-adic numbers, Characterization
theorem, Skitovich-Darmois theorem

\bigskip

\textbf{Mathematics Subject Classification}: 60B15 · 62E10 · 43A35

\bigskip

\textbf{1. Introduction}

\bigskip

One of the most important characterization theorem of mathematical
statistics is the Kac--Bernstein theorem. This theorem characterizes
a Gaussian distribution by the independence of the sum and of the
difference of two independent random variables. V.P. Skitovich and
G. Darmois generalized this theorem.

\textbf{Theorem A} (\cite{Skitovich}, \cite{Darmois}, see also
\cite[Ch. 3]{Kag-Lin-Rao}). \textsl{Let $\xi_j, \ j=1, 2,\dots, n,\
n\geq 2$, be independent random variables. Let  $\alpha_j, \beta_j$
be nonzero constants. If the linear forms
$L_1=\alpha_1\xi_1+\cdots+\alpha_n\xi_n$ and
$L_2=\beta_1\xi_1+\cdots+\beta_n\xi_n$ are independent, then all
random variables $\xi_j$ are Gaussian.}

The classical characterization theorems of mathematical statistics
were extended to different algebraic structures such as
finite-dimensional and infinite-dimensional linear spaces, symmetric
spaces, quantum groups (see e.g. \cite{Mi2},
\cite{Neue}--\cite{NeueScho}). Much attention has been also devoted
to the study of analogues of the Skitovich--Darmois theorem for some
classes of locally compact Abelian groups. In this case coefficients
of linear forms are topological automorphisms of a group. (see e.g.
\cite{Fe10}--\cite{Fe14},
\cite{Fe2013}--\cite{FeGra2},\cite{Mazur1}-- \cite{FelMyr2014},
\cite{Stapleton} and also \cite{Fe0}, where one can find additional
references). Note that most articles were devoted to the description
of groups on which the independence of two linear forms from $n$
random variables implies that all random variables are Gaussian or
have distributions which can be considered as analogues of Gaussian
distributions on groups (e.g. idempotent distributions). We also
note that the description of distributions on groups which are
characterized by the independence of linear forms depends not only
on a group but also on the length of forms. For example, on the one
hand, it was proved that for $n=2$ an analogue of the
Skitovich--Darmois theorem is valid for all finite groups
(\cite{Fe14}). On the other hand, if $n>2$ then the analogue of the
Skitovich--Darmois theorem is valid only for some special finite
Abelian groups (\cite{Fe10}). Similar results were obtained for
random variables with values in discrete Abelian groups (see
\cite{FeGra2} and \cite{Fe13}). The case of compact Abelian groups
was studied in \cite{Fe12} and \cite{FeGra1}.

The situation is essentially changing, if we consider  $n$ linear
forms of $n$ independent random variables instead of two linear
forms of $n$ independent random variables. We note that one of the
first paper in which  $n$ linear forms were considered, was the
paper \cite{Stapleton}, where J.H.Stapleton considered the case of
compact connected Abelian groups with some restrictions on
coefficients of the forms and random variables. Next it was proved
in \cite{Mazur1}, that if random variables take values in a finite
Abelian group then the independence of $n$ linear forms of $n$
independent random variables implies that random variables are
idempotent. The cases of discrete and compact Abelian groups were
studied in \cite{Mazur2} and \cite{Mazur3}.

The distributions on a group of $p$-adic numbers $\Omega_p$, which
are characterized by the independence of two linear forms of two
independent random variables were studied in \cite{Fe2013}. Since
the group $\Omega_p$ is totally disconnected, the Gaussian
distributions on it are degenerated. The role of Gaussian
distributions on $\Omega_p$ plays idempotent distributions, i.e. the
set of shifts of the Haar distributions of the compact subgroups. It
turned out that on this group the description of distributions which
are characterized by the independence of two linear forms of two
independent random variables depends on the coefficients of the
linear forms. In particular, all the coefficients of linear forms
for which the independence of two linear forms of two independent
random variables implies that random variables are idempotent were
described in \cite{Fe2013}.

In this paper we study distributions on the group of $p$-adic
numbers $\Omega_p$, which are characterized by the independence of
three linear forms of three independent random variables. This
studying is much more complicated then in \cite{Fe2013} because the
number of parameters increases. The main result of the paper is
Theorem 1. In particular, we obtain in Theorem 1 the description of
all the coefficients of linear forms for which the independence of
three linear forms of three independent random variables implies
that all random variables are idempotent.

\bigskip

\textbf{2. Definitions and notation}

\bigskip

We use some results of the duality theory for the locally compact
Abelian groups (see  \cite{HeRo}). For an arbitrary locally compact
Abelian group $X$  let $Y=X^\ast$ be its character group, and
$(x,y)$  be the value of a character $y \in Y$ at an element $x \in
X$. If $K$ is a closed subgroup of $X$, we denote by $A(Y, K) = \{y
\in Y: (x, y) = 1$ for all $x \in K \}$ its annihilator. If $\delta
: X \mapsto X$ is a continuous homomorphism, then the adjoint
homomorphism $\widetilde \delta : Y \mapsto Y$ is defined by the
formula $(x, \widetilde \delta y) = (\delta x, y)$ for all $x \in X,
\ y \in Y$. Denote by $Aut(X)$ the group of topological
automorphisms of the group $X$. We note that $\delta \in {\rm
Aut}(X)$ if and only if $\widetilde\delta \in {\rm Aut}(Y)$. Denote
by $I$ the identity automorphism of a group. Denote by $f_n$
homomorphism defined by the formula $f_n(x)=nx$. Put $X_{(n)}=\{x\in
X: nx=0\}$.

Let  ${M^1}(X)$ be the convolution semigroup of probability
distributions on $X.$ For a distribution $\mu \in {M^1}(X)$ denote
by $$\widehat \mu(y) = \int_X (x, y) d \mu(x)$$ its characteristic
function (Fourier transform), and by $\sigma(\mu)$ the support of
$\mu$. For $\mu \in {M^1}(X)$, we define the distribution $\bar \mu
\in M^1(X)$ by the formula $\bar \mu(E) = \mu(-E)$ for any Borel set
$E \subset X$. Observe that $\widehat {\bar \mu}(y) =
\overline{\widehat \mu(y)}$. Let $K$ be a compact subgroup of  $X$.
Denote by $m_K$  the Haar distribution on $K$. We note that the
characteristic function of $m_K$ is of the form
\begin{equation}\label{i2}
    \widehat m_K(y) = \cases{ 1, \quad  \quad y \in A(Y,K), \cr\cr 0,
    \quad  \quad y \notin A(Y,K).  \cr}
\end{equation}
Denote by $I(X)$ the set of all idempotent distributions on $X$,
i.e. the set of shifts of the Haar distributions $m_K$ of the
compact subgroups $K$ of $X$. Let $x\in X$. Denote by $E_x$ the
degenerated distribution concentrated at the point $x$. Denote by
$D(X)$ the set of all degenerated distributions on $X$.

Let $p$ be a prime number. We need some properties of the group of
$p$-adic numbers $\Omega_p$ (see  \cite[\S 10]{HeRo}). As a set
$\Omega_p$ coincides with the set of sequences of integers of the
form $x=(\dots,x_{-n}, x_{-n+1},\dots, x_{-1}, x_0, x_1,\dots,x_n,
x_{n+1},\dots),$ where $x_n \in\{0, 1,\dots, p-1\}$, such that
$x_n=0$ for $n < n_0$, where the number $n_0$ depends on $x$. We
correspond to each element $x \in \Omega_p$ the series
$\sum\limits_{k=-\infty}^{\infty} x_k p^k.$ Addition and
multiplication of the series are defined in a natural  way and they
define the operations of addition and  multiplication in $\Omega_p$.
With respect to these operations $\Omega_p$ is a field. Denote by
$\Delta_p$ a subgroup of $\Omega_p$ consisting of $x \in \Omega_p$
such that $x_n=0$ for $n < 0$. The subgroup $\Delta_p$ is called the
group of $p$-adic integers. The family of the subgroups
$\{p^k\Delta_p\}_{k=-\infty}^{\infty}$ forms an open basis at zero
of the group $\Omega_p$ and defines a topology on $\Omega_p$. With
respect to this topology the group $\Omega_p$ is locally compact,
non-compact, and totally disconnected. We note that the group
$\Omega_p$ is represented as a union $\Omega_p=
\bigcup\limits_{k=-\infty}^{\infty} p^k \Delta_p$. Note that any
compact subgroup $K$ of $\Omega_p$ is of the form $K=p^l \Delta_p$
for some $l$. The character group $\Omega_p^{*}$ of the group
$\Omega_p$ is topologically isomorphic to $\Omega_p$. Each
automorphism $\alpha \in {\rm Aut}(\Omega_p)$ is of the form $\alpha
g=x_{\alpha} g,$ $g \in \Omega_p,$ where $x_\alpha \in \Omega_p,$
$x_\alpha \ne 0$. For $\alpha \in {\rm Aut}(\Omega_p)$ we  identify
the automorphism $\alpha \in {\rm Aut}(\Omega_p)$ with the
corresponding element $x_\alpha \in \Omega_p$, i.e. when we write
$\alpha g$, we  suppose that $\alpha \in \Omega_p$. We note that
$\widetilde\alpha=\alpha$. Denote by $\Delta_p^0$ the subset of
$\Omega_p$  consisting of all invertible in $\Delta_p$ elements,
$\Delta_p^0=\{x \in \Omega_p: x_n =0$ for $n < 0, \ x_0 \ne 0 \}$.
We note that each element  $g \in \Omega_p$ is represented in the
form $g = p^k c,$ where $k$ is an integer, and $c \in \Delta_p^0$.
Hence, multiplication on $c$ is a topological automorphism of the
group $\Delta_p$.

Denote by ${\mathbb Z}(p^\infty)$ the set of rational numbers of the
form $\{{k / p^n} : k=0, 1, \dots,p^n-1, \ n=0,1,\dots\}$. If we
define the operation in ${\mathbb Z}(p^\infty)$ as addition modulo
1, then ${\mathbb Z}(p^\infty)$ is transformed into an Abelian group
which we consider in the discrete topology. Obviously, this group is
topologically isomorphic to the multiplicative group of all $p^n$-th
roots of unity, where $n$ goes through the set of nonnegative
integers considering in the discrete topology. For a fixed $n$
denote by ${\mathbb Z}(p^n)$ a subgroup of ${\mathbb Z}(p^\infty)$
consisting of all elements of the form ${\{{k / p^n} : k=0, 1,
\dots,p^n-1\}}$. Note that the group ${\mathbb Z}(p^n)$ is
topologically isomorphic to the multiplicative group of all $p^n$th
roots of unity considering in the discrete topology. Observe that
the groups ${\mathbb Z}(p^\infty)$ and $\Delta_p$ are the character
groups of one another.

\bigskip

\textbf{3. Main results}

\bigskip

Let $X = \Omega_p$, and let  $\xi_1, \xi_2, \xi_3$ be independent
random variables with values in $X$ and distributions $\mu_1, \mu_2,
\mu_3$. Let $\alpha_j, \beta_j, \gamma_j \in {\rm Aut}(X)$. We
consider linear forms $L_1 = \alpha_1\xi_1 + \alpha_2 \xi_2+
\alpha_3 \xi_3$, $L_2=\beta_1\xi_1 + \beta_2 \xi_2+ \beta_3 \xi_3$
and $L_3=\gamma_1\xi_1 + \gamma_2 \xi_2+ \gamma_3 \xi_3$ and assume
that $L_1$, $L_2$ and $L_3$ are independent.

We can consider new independent random variables $\xi'_1 =
\alpha_1\xi_1$, $\xi'_2 = \alpha_2\xi_2$, $\xi'_3 = \alpha_3\xi_3$,
and reduce the problem of description of possible distributions
$\mu_1, \mu_2, \mu_3$ to the case when
$$L_1 = \xi_1 + \xi_2 + \xi_3,$$
$$L_2=\alpha'_1\xi_1 + \alpha'_2 \xi_2+ \alpha'_3 \xi_3,$$
$$L_3=\beta'_1\xi_1 + \beta'_2 \xi_2+ \beta'_3 \xi_3.$$
Since $\alpha'_j, \beta'_j \in {\rm Aut}(X)$, we have
$\alpha'_j=p^{k_j} c_j, \beta'_j=p^{l_j} d_j$, where all $c_j, d_j
\in \Delta_p^0$, $k_j, l_j \in \mathbb{Z}$.

By renumbering random variables, we can assume that there are two
possible cases:

1. $k_1=min\{k_1, k_2, k_3\}$, $l_1=min\{l_1, l_2, l_3\}$.

2. $k_1=min\{k_1, k_2, k_3\}$, $l_2=min\{l_1, l_2, l_3\}$.

Since $L_1$, $L_2$ and $L_3$ are independent if and only if either
$L_1$, $L'_2=\alpha_1^{-1} L_2$ and $L'_3=\beta_1^{-1} L_3$, or
$L_1$, $L'_2=\alpha_1^{-1} L_2$ and $L'_3=\beta_2^{-1}L_3$ are
independent, the problem of description of possible distributions
$\mu_1, \mu_2, \mu_3$ is reduced to the case when
$$L_1 = \xi_1 + \xi_2 + \xi_3,$$
$$L_2=\xi_1 + \delta_1 \xi_2+ \delta_2 \xi_3,$$
$$L_3=\xi_1 + \varepsilon_1 \xi_2+ \varepsilon_2 \xi_3,$$
or to the case when
$$L_1 = \xi_1 + \xi_2 + \xi_3,$$
$$L_2=\delta_1\xi_1 +  \xi_2+ \delta_2 \xi_3,$$
$$L_3=\xi_1 + \varepsilon_1\xi_2+ \varepsilon_2 \xi_3.$$
Note that it follows from the choice of automorphisms that
\begin{equation}\label{t1}
    \delta_j=p^{k_j} c_j, \varepsilon_j=p^{l_j} d_j,
    \quad c_j, d_j \in \Delta_p^0, k_j\geq 0, l_j\geq 0, j=1, 2,
\end{equation}
in both cases. Let $k=min \{ k_1,k_2,l_1,l_2\}$.

Denote
$$\Lambda_1= \left(%
\begin{array}{ccc}
  1 & 1 & 1 \\
  1 & \delta_1 & \delta_2 \\
  1 & \varepsilon_1 & \varepsilon_2 \\
\end{array}%
\right), \quad \Lambda_2= \left(%
\begin{array}{ccc}
  1 & 1 & 1 \\
  \delta_1 & 1 & \delta_{2} \\
  1 & \varepsilon_1 & \varepsilon_2 \\
\end{array}
\right)$$

Let $\Lambda=\Lambda_1$ or $\Lambda=\Lambda_2$. We will keep the
designation $\Lambda_1$ and $\Lambda_2$ throughout the article. Then
$\det \Lambda= p^q \lambda$, where $q\geq 0$, $\lambda \in
\Delta_p^0$.

The main result of the article is to get necessary and sufficient
conditions on coefficients of the linear forms $L_j$ when the
independence of $L_j$ implies that all $\mu_j$ are independent.

\medskip

\textbf{Theorem 1.}
\textsl{Let $X = \Omega_p$. Let $\xi_1, \xi_2, \xi_3$ be independent
random variables with values in $X$ and distributions $\mu_1, \mu_2,
\mu_3$. Let $\delta_j, \varepsilon_j \in {\rm Aut}(X)$, where
$\delta_j=p^{k_j} c_j, \varepsilon_j=p^{l_j} d_j$, where $c_j, d_j
\in \Delta_p^0$, $ k_j\geq 0, l_j\geq 0, j=1, 2$. Put $k=min \{
k_1,k_2,l_1,l_2\}$.}

\textsl{Consider a random vector $(L_1,L_2,L_3)= \Lambda
(\xi_1,\xi_2,\xi_3)$, where either $\Lambda=\Lambda_1$ or
$\Lambda=\Lambda_2$. Let $\det \Lambda=p^q \lambda$, where $q\geq
0$, $\lambda \in \Delta_p^0$. Suppose that the random vector
$(L_1,L_2,L_3)$ has independent components. Then the following
statements hold.}

\textsl{\textbf{A.} Assume that $q=0$.}

\textsl{All $\mu_j\in I(X)$ iff each row and each column of the
matrix $\Lambda$ contains at least two elements from $\Delta_p^0$.}

\textsl{\textbf{B.} Assume that $q>0$.}

\textsl{\textbf{B.1.} Assume that at least one of the elements
$\delta_j, \varepsilon_j\in \Delta_p^0$. Then all $\mu_j\in D(X)$.}

\textsl{\textbf{B.2.} Assume that all $\delta_j,
\varepsilon_j\not\in \Delta_p^0$. All $\mu_j\in I(X)$ iff either
$q>k$ (in this case all $\mu_j\in D(X)$) or $q=k$ and either
$k_1=l_2=k$ or $k_2=l_1=k$.}

\medskip

\textbf{Remark 1.} Case B.2 is possible only when
$\Lambda=\Lambda_1$.

\medskip

To prove Theorem 1 we need some lemmas.

\medskip

\textbf{Lemma 1} (\cite{Mazur1}).
\textsl{Let $X$ be a second countable locally compact Abelian group,
$Y$ be its character group. Let $\alpha_j, \beta_j, \gamma_j$ be
continuous endomorphisms of $X$, $\xi_1, \xi_2, \xi_3$ be
independent random variables with values in $X$ and distributions
$\mu_1, \mu_2, \mu_3$. Let $L_1 = \alpha_1\xi_1 + \alpha_2 \xi_2+
\alpha_3 \xi_3$, $L_2=\beta_1\xi_1 + \beta_2 \xi_2+ \beta_3 \xi_3$
and $L_3=\gamma_1\xi_1 + \gamma_2 \xi_2+ \gamma_3 \xi_3$. The linear
forms $L_1$, $L_2$ and $L_3$ are independent iff the characteristic
functions $\hat\mu_j(y)$ satisfy the equation}
\begin{equation}\label{1}
    \hat\mu_1(\tilde{\alpha_1}u +  \tilde{\beta_1}v + \tilde{\gamma_1} w) \hat\mu_2(\tilde{\alpha_2}u +  \tilde{\beta_2}v + \tilde{\gamma_2} w)
    \hat\mu_3(\tilde{\alpha_3}u +  \tilde{\beta_3}v + \tilde{\gamma_3} w)=$$ $$=
    \hat\mu_1(\tilde{\alpha_1} u)\hat\mu_2(\tilde{\alpha_2} u)\hat\mu_3(\tilde{\alpha_3} u)
    \hat\mu_1(\tilde{\beta_1} v)\hat\mu_2(\tilde{\beta_2} v)\hat\mu_3(\tilde{\beta_3} v)
    \hat\mu_1(\tilde{\gamma_1} w)\hat\mu_2(\tilde{\gamma_2} w)\hat\mu_3(\tilde{\gamma_3} w), \quad u, v, w \in
    Y.
\end{equation}
\medskip

The proof of the next lemma is analogous to the proof of Lemma 2 of
the paper \cite{Fe2013}.

\medskip

\textbf{Lemma 2.}
\textsl{Let $X = \Omega_p$, $Y$ be its character group. Let $\xi_1,
\xi_2, \xi_3$ be independent random variables with  values in $X$
and distributions $\mu_1, \mu_2, \mu_3$ such that $\hat\mu_j(y) \ge
0$, $j = 1, 2, 3$. Let $\alpha_j, \beta_j \in {\rm Aut}(X)$ and
$\alpha_j=p^{k_j} c_j, \beta_j=p^{l_j} d_j \in {\rm Aut}(X)$, where
all $k_j\geq 0$, $l_j\geq 0$ and $c_j, d_j\in \Delta_p^0$. Assume
that the linear forms $L_1 = \xi_1 + \xi_2 + \xi_3$,
$L_2=\alpha_1\xi_1 + \alpha_2 \xi_2+ \alpha_3 \xi_3$ and
$L_3=\beta_1\xi_1 + \beta_2 \xi_2+ \beta_3 \xi_3$ are independent.
Then there exists a subgroup $B= p^n\Delta_p$ of the group $Y$ such
that all $\sigma(\mu_j) \subset B$.}

\textbf{Proof.} We use the fact that the family of the subgroups
$\{p^l \Delta_p\}_{l=-\infty}^\infty$ forms an open basis at zero of
the group $Y$. Since
$\widehat\mu_1(0)=\widehat\mu_2(0)=\widehat\mu_3(0)=1$, we can
choose $m\geq 0$ in such a way that  $\widehat\mu_j(y) >0$ for $y
\in L=p^m \Delta_p,$ $j=1, 2, 3$. Note that $\alpha_j(L) \subset L$,
$j=1, 2, 3$. Put $\psi_j(y)= - \log\widehat\mu_j(y),$ $y \in L,$
$j=1, 2, 3$.

By Lemma 1 the characteristic functions $\widehat\mu_j(y)$ satisfy
equation (\ref{1}). Taking into account that
$\widetilde\alpha_j=\alpha_j$, $\widetilde\beta_j=\beta_j$, we get
from (\ref{1}) the following equation
\begin{equation}\label{l2.1}
    \psi_1(u+v+w)+\psi_2(\alpha_1 u+\alpha_2 v+\alpha_3 w)+ \psi_3(\beta_1 u+\beta_2 v+\beta_3 w)=$$
    $$= \psi_1(u)+\psi_2(\alpha_1 u)+\psi_3(\beta_1 u) +\psi_1(v)+\psi_2(\alpha_2 v)+\psi_3(\beta_2 v)
    +\psi_1(w)+\psi_2(\alpha_3 w)+\psi_3(\beta_3 w), \quad u,v,w \in L.
\end{equation}
Put $w=0$ in (\ref{l2.1}) and integrate the obtained equation over
the group $L$ with respect to the Haar distribution $dm_L(u)$. Using
the fact that the Haar distribution $m_L$ is $L$-invariant, we
obtain
\begin{equation}\label{l2.2}
    \psi_1(v)+\psi_2(\alpha_2 v)+\psi_3(\beta_2 v)=0, \quad v \in L.
\end{equation}
Since $\psi_j(y)\geq 0$ on $Y$, it follows from (\ref{l2.2}) that
$\psi_1(v)=\psi_2(\alpha_2 v)=\psi_3(\beta_2 v)=0$ for $v \in L$.
Thus $\widehat\mu_1(y)=\widehat\mu_2(\alpha_2
y)=\widehat\mu_3(\beta_2 y)=1, \ y \in L$. Put $B=L\cap \alpha_2(L)
\cap \beta_2(L)$, i.e. $B= p^n\Delta_p$, where $n=max\{k_2, l_2\}$.
Then $B$ is a required subgroup. Lemma 2 is proved.

$\blacksquare$

\bigskip
\textbf{Lemma 3 (\cite{Mazur2}).}
\textsl{Let $X = \Delta_p$, let $\xi_i$, $i=1,...,n$, be independent
random variables with values in $X$ and with distributions $\mu_i$.
Let $\alpha_{ij} \in {\rm Aut}(X)$. If the linear form
$L_j=\sum_{1}^{n} \alpha_{ij} \xi_i$, $j=1,...,n$, are independent
then all $\mu_i\in I(X)$.}

\bigskip

The following lemma is a partial case of the general theorem of the
paper \cite{Mazur3}.

\bigskip
\textbf{Lemma 4 (\cite{Mazur3}).}
\textsl{Let $X$ be a finite Abelian group. Let $\alpha_j, \beta_j,
\gamma_j$ be endomorphisms of the group $X$. Put
\begin{equation}\label{3}
    \Lambda=\left(%
            \begin{array}{ccc}
              \alpha_1 & \alpha_2 & \alpha_3 \\
               \beta_1 & \beta_2 & \beta_3 \\
              \gamma_1 & \gamma_2 & \gamma_3 \\
            \end{array}
\right).
\end{equation}
Assume that each column of the matrix $\Lambda$ contain at least two
automorphisms of the group $X$, and the matrix $\Lambda\in
Aut(X^3)$.}

\textsl{Let $\xi_1, \xi_2, \xi_3$ be independent random variables
with values in $X$ and distributions $\mu_1, \mu_2, \mu_3$ such that
each support $\sigma(\mu_j)$ is not contained in a coset of a proper
subgroup of the group $X$. If the linear forms $L_1 = \alpha_1\xi_1
+ \alpha_2 \xi_2+ \alpha_3 \xi_3$, $L_2=\beta_1\xi_1 + \beta_2
\xi_2+ \beta_3 \xi_3$ and $L_3=\gamma_1\xi_1 + \gamma_2 \xi_2+
\gamma_3 \xi_3$ are independent, then all distributions
$\mu_1=\mu_2=\mu_3=m_X$.}

\bigskip

It is convenient for us to formulate as a lemma the following
well-known statement (see e.g. \cite[Proposition 2.13]{Fe0}).

\bigskip
\textbf{Lemma 5.} \textsl{Let $X$ be a locally compact Abelian
group, $Y$ be its character group. Let $\mu\in{\rm M}^1(X)$. Then
the set $E=\{y\in Y:\ \hat\mu(y)=1\}$ is a closed subgroup of $Y$,
the characteristic function $\hat\mu(y)$ is $E$-invariant, i.e.,
$\hat\mu(y)$ takes a constant value on each coset of the group $Y$
with respect to $E$, and $\sigma(\mu)\subset A(X,E)$.}

\bigskip
\textbf{Lemma 6.}
\textsl{Let $X = \mathbb{Z}(p^n)$, $Y$ be its character group. Let
$\alpha_j, \beta_j, \gamma_j$ be endomorphisms of the group $X$. Put
\begin{equation}\label{3}
    \Lambda=\left(%
            \begin{array}{ccc}
               \alpha_1 & \alpha_2 & \alpha_3 \\
               \beta_1 & \beta_2 & \beta_3 \\
              \gamma_1 & \gamma_2 & \gamma_3 \\
            \end{array}
\right).
\end{equation}
Assume that each row and column of the matrix $\Lambda$ contains at
least two automorphisms of the group $X$, and the matrix $\Lambda\in
Aut(X^3)$.}

\textsl{Let $\xi_1, \xi_2, \xi_3$ be independent random variables
with values in $X$ and distributions $\mu_1, \mu_2, \mu_3$ such that
$\hat\mu_j(y) \ge 0$, $j = 1, 2, 3$, and at least one support
$\sigma(\mu_j)$ is not contained in a coset of a proper subgroup of
$X$. If the linear forms $L_1 = \alpha_1\xi_1 + \alpha_2 \xi_2+
\alpha_3 \xi_3$, $L_2=\beta_1\xi_1 + \beta_2 \xi_2+ \beta_3 \xi_3$
and $L_3=\gamma_1\xi_1 + \gamma_2 \xi_2+ \gamma_3 \xi_3$  are
independent, then $\mu_1=\mu_2=\mu_3=m_X$.}

\bigskip

\textbf{Proof.} Since $X = \mathbb{Z}(p^n)$, we have $Y \approx
\mathbb{Z}(p^n)$. To avoid introducing new notation we will suppose
that $Y = \mathbb{Z}(p^n)$. Lemma 1 implies that the characteristic
functions $\hat\mu_j(y)$ satisfy equation (\ref{1}).

Suppose that there exists such element $y_0\in Y$ that one of the
functions $\hat\mu_j(y)$ is equal to 1 at this element. We may
assume without loss of generality that $\hat\mu_1(y_0)=1$. Since the
set of elements where a characteristic function is equal to one is a
subgroup, and any nonzero subgroup of $Y$ contains the subgroup
$\mathbb{Z}(p)$, we obtain that $\hat\mu_1(y)=1$ for
$y\in\mathbb{Z}(p)$. Since $\mathbb{Z}(p)$ is a characteristic
subgroup, we can consider the restriction of equation (\ref{1}) to
$\mathbb{Z}(p)$. Thus, we have

\begin{equation}\label{l4.1}
    \hat\mu_2({\alpha_2}u +  {\beta_2}v + {\gamma_2} w)
    \hat\mu_3({\alpha_3}u +  {\beta_3}v + {\gamma_3} w)=$$ $$=
    \hat\mu_2({\alpha_2} u)\hat\mu_3({\alpha_3} u)
    \hat\mu_2({\beta_2} v)\hat\mu_3({\beta_3} v)
    \hat\mu_2({\gamma_2} w)\hat\mu_3({\gamma_3} w), \quad u, v, w \in
    \mathbb{Z}(p).
\end{equation}
Suppose that all coefficients $\alpha_j, \beta_j, \gamma_j$
($j=2,3$) are automorphisms of the group $Y$. It is obvious that
there exists a non-zero element $(u_0,v_0,w_0)$ such that
\begin{equation}\label{l4.2}
    \left\{%
\begin{array}{ll}
    {\alpha_2}u_0 +  {\beta_2}v_0 + {\gamma_2} w_0=0 & \hbox{;} \\
    {\alpha_3}u_0 +  {\beta_3}v_0 + {\gamma_3} w_0=0 & \hbox{.} \\
\end{array}%
\right.
\end{equation}
Put $u=u_0$, $v=v_0$, $w=w_0$ in (\ref{l4.1}). We obtain
\begin{equation}\label{l4.3}
    1=
    \hat\mu_2({\alpha_2} u_0)\hat\mu_3({\alpha_3} u_0)
    \hat\mu_2({\beta_2} v_0)\hat\mu_3({\beta_3} v_0)
    \hat\mu_2({\gamma_2} w_0)\hat\mu_3({\gamma_3} w_0).
\end{equation}
Since the element $(u_0,v_0,w_0)$ is non-zero, and the set of
elements where a characteristic function is equal to 1 is a
subgroup, it follows from (\ref{l4.3}) that
$\hat\mu_2(y)=\hat\mu_3(y)=1$ for $y\in\mathbb{Z}(p)$. Thus Lemma 5
implies that all supports $\sigma(\mu_j)\subset A(X,\mathbb{Z}(p))$.
We obtain the contradiction to the conditions of the lemma.

Suppose that the only one of the coefficients $\alpha_j, \beta_j,
\gamma_j$ ($j=2,3$) is not an automorphism, and all the rest are
automorphisms of the group $Y$. Without loss of generality we can
assume that $\gamma_3\not\in Aut(X)$. Then equation (\ref{l4.1})
takes the form
\begin{equation}\label{l4.4}
    \hat\mu_2({\alpha_2}u +  {\beta_2}v + {\gamma_2} w)
    \hat\mu_3({\alpha_3}u +  {\beta_3}v )=
    \hat\mu_2({\alpha_2} u)\hat\mu_3({\alpha_3} u)
    \hat\mu_2({\beta_2} v)\hat\mu_3({\beta_3} v)
    \hat\mu_2({\gamma_2} w), \quad u, v, w \in
    \mathbb{Z}(p).
\end{equation}
Put $u=-\alpha_3^{-1}\beta_3 y$, $v=y$,
$w=\gamma_2^{-1}(\alpha_2\alpha_3^{-1}\beta_3-\beta_2)y$ in
(\ref{l4.4}). We obtain
\begin{equation}\label{l4.5}
    1=
    \hat\mu_2(-\alpha_3^{-1}\beta_3{\alpha_2} y)\hat\mu_3(-\beta_3 y)
    \hat\mu_2({\beta_2} y)\hat\mu_3({\beta_3} y)
    \hat\mu_2((\alpha_2\alpha_3^{-1}\beta_3-\beta_2)y), \quad y \in
    \mathbb{Z}(p).
\end{equation}
It follows from this that $\hat\mu_2(y)=\hat\mu_3(y)=1$ for
$y\in\mathbb{Z}(p)$. Thus Lemma 5 implies that all supports
$\sigma(\mu_j)\subset A(X,\mathbb{Z}(p))$. We obtain the
contradiction to the conditions of the lemma.

Suppose that two of the coefficients $\alpha_j, \beta_j, \gamma_j$
($j=2,3$) are not automorphisms, and all the rest coefficients are
automorphisms of the group $Y$. By the conditions on $\Lambda$ it is
possible when automorphisms are coefficients for different
variables. Without loss of generality we can assume that $\beta_2,
\gamma_3\not\in Aut(X)$. Then equation (\ref{l4.1}) takes the form
\begin{equation}\label{l4.6}
    \hat\mu_2({\alpha_2}u +  {\gamma_2} w)
    \hat\mu_3({\alpha_3}u +  {\beta_3}v )=
    \hat\mu_2({\alpha_2} u)\hat\mu_3({\alpha_3} u)
    \hat\mu_3({\beta_3} v)
    \hat\mu_2({\gamma_2} w), \quad u, v, w \in
    \mathbb{Z}(p).
\end{equation}
Put $u=y$, $v=-{\beta_3}^{-1}{\alpha_3}y$, $w=
-{\gamma_2}^{-1}{\alpha_2}y$ in  (\ref{l4.6}). We obtain
\begin{equation}\label{l4.7}
    1=
    \hat\mu_2({\alpha_2} y)\hat\mu_3({\alpha_3} y)
    \hat\mu_3(-{\alpha_3}y)
    \hat\mu_2(-{\alpha_2}y), \quad y \in
    \mathbb{Z}(p).
\end{equation}
It follows from this that $\hat\mu_2(y)=\hat\mu_3(y)=1$ for
$y\in\mathbb{Z}(p)$. Thus, Lemma 5 implies that all supports
$\sigma(\mu_j)\subset A(X,\mathbb{Z}(p))$. We obtain the
contradiction to the conditions of the lemma.

We have considered all possibilities.

Thus we obtain that all $\hat\mu_j(y)<1$ for $y\in Y\setminus
\{0\}$. Now we are in the conditions of Lemma 4. The statement of
the lemma follows from Lemma 4.

$\blacksquare$

\bigskip

\textbf{Lemma 7.} \textsl{Let $X$ be a locally compact Abelian
group, $Y$ be its character group. Let $H$ be an open subgroup of
$Y$. Let $f_0(y)$ be a continuous positive definite function on $H$,
and $f(y$ be the function on $Y$ such that}
\begin{equation}\label{l8.1}
    f(y)= \left\{%
\begin{array}{ll}
    f_0(y), & \hbox{$y\in H$;} \\
    0, & \hbox{$y\not\in H$.} \\
\end{array}%
\right.
\end{equation}
\textsl{Then $f(y)$ is a continuous positive definite function, and
there exists a distribution $\mu\in M^1(X)$ such that
$\hat\mu(y)=f(y)$.}

\bigskip

\textbf{Proof.} The function $f(y)$ is a positive definite function
on $Y$ (\cite[\S 32]{HeRo2}]). Since $H$ is open, the function
$f(y)$ is continuous/ By the Bochner theorem there exists a
distribution $\mu$ such that $\widehat\mu(y) = f(y)$.

$\blacksquare$

\bigskip

The following lemma is obvious and follows immediately from the form
of endomorphisms of $\Delta_p$. We will constantly use this lemma,
but not refer to it.

\bigskip

\textbf{Lemma 8.} \textsl{Let $X = \Delta_p$. Let $\alpha\in
Aut(X)$, let $\beta$ be an endomorphism of $X$ such that
$\beta\not\in Aut(X)$. Then $\alpha-\beta\in Aut(X)$.}

\bigskip



\textbf{Lemma 9.} \textsl{Let $X = \Delta_p$. Let $\xi_1, \xi_2,
\xi_3$ be independent random variables with values in $X$ and with
distributions $\mu_1, \mu_2, \mu_3$ such that $\hat\mu_j(y)\geq 0$,
$j=1,2,3$. Suppose that at least one of the supports $\sigma(\mu_j)$
is not contained in a proper subgroup of the group $X$. Let $c_j,
d_j \in {\rm Aut}(X)$. Consider linear forms $L_1 = \xi_1 + \xi_2+
\xi_3$, $L_2 = \xi_1 + p^{k_1} c_1 \xi_2+ p^{k_2} c_2 \xi_3$, and
$L_3=\xi_1 + p^{l_1} d_1 \xi_2+ p^{l_2} d_2 \xi_3$, where all
$k_j>0, l_j>0$. Put $k=min \{ k_1,k_2,l_1,l_2\}$. Suppose that
$L_1$, $L_2$ and $L_3$ are independent. All $\mu_j\in I(X)$ iff
either $k_1=l_2=k$ or $k_2=l_1=k$. }

\bigskip

\textbf{Proof.}

Let $l, m, n \in \mathbb{N}$ such that $0<k<l<m<n$. Up to notation
of the random variables $\xi_2$ and $\xi_3$ and up to permutation of
the linear forms $L_2$ and $L_3$ there are possible the following
cases for $k_j, l_j$.

\small

\noindent\begin{tabular}{|c|c|c|c|c|c|c|c|c|c|c|c|c|c|c|c|c|c|c|c|c|c|}
  \hline
  . & 1 & 2 & 3 & 4 & 5 & 6 & 7 & 8 & 9 & 10 & 11 & 12 & 13 & 14 & 15 & 16 & 17 & 18 & 19 & 20 & 21 \\
  \hline
  $k_1$ & k & k & k & k & k & k & k & k & k & k & k & k & k & k & k & k & k & k & k & k & k \\
  $l_1$ & k & k & k & l & l & k & l & l & l & l & l & m & l & m & m & l & l & m & m & n & n \\
  $k_2$ & k & k & l & k & l & l & k & m & l & l & m & l & m & l & m & m & n & l & n & l & m \\
  $l_2$ & k & l & l & l & k & m & m & k & l & m & l & l & m & m & l & n & m & n & l & m & l \\
  \hline
\end{tabular}

\normalsize

Since $X =\Delta_p$, we have $Y \approx \mathbb{Z}(p^{\infty})$. To
avoid introducing new notation we will suppose that $Y =
\mathbb{Z}(p^{\infty})$.

Put $f(y)=\hat\mu_1(y)$, $g(y)=\hat\mu_2(y)$, $h(y)=\hat\mu_3(y)$.
By Lemma 1 the functions $f(y)$, $g(y)$, $h(y)$ satisfy (\ref{1})
which takes the form

\begin{equation}\label{l5.1}
f(u+v+w) g(u+p^{k_1}c_1v+p^{l_1}d_1w)h(u+p^{k_2}c_2v+p^{l_2}d_2w)=$$
$$=f(u)f(v)f(w)
g(u)g(p^{k_1}c_1v)g(p^{l_1}d_1w)h(u)h(p^{k_1}c_2v)h(p^{l_2}d_2w),
u,v,w\in Y,
\end{equation}
where $f(y)=\hat\mu_1(y)$, $g(y)=\hat\mu_2(y)$, $h(y)=\hat\mu_3(y)$,
holds.

We note that if at least one of the supports $\sigma(\mu_j)$ is not
contained in a proper subgroup of $X$, then by Lemma 5 $\{y\in Y:
f(y)=g(y)=h(y)=1\}=\{0\}$.

Put
$$\Lambda= \left(%
\begin{array}{ccc}
  1 & 1 & 1 \\
  1 & p^{k_1}c_1 & p^{k_2}c_2 \\
  1 & p^{l_1}d_1 & p^{l_2}d_2 \\
\end{array}%
\right).$$ We have
\begin{equation}\label{l5.6}
    \det
\Lambda=p^{k_1+l_2}c_1d_2-p^{k_2+l_1}c_2d_1+p^{k_2}c_2-p^{k_1}
c_1+p^{l_1}d_1-p^{l_2}d_2.
\end{equation}
Note that $\det\Lambda$ can be represented in the form
\begin{equation}\label{l5.7}
    \det \Lambda= p^q \lambda,
\end{equation}
where $\lambda\in Aut(Y)$, $q>0$. In this case $q\geq k$.

Putting $u=(p^{k_1+l_2}c_1d_2-p^{k_2+l_1}c_2d_1)y$,
$v=(p^{l_2}d_2-p^{l_1}d_1)y$, $w=(p^{k_2}c_2-p^{k_1}c_1)y$, $y\in
Y$, in (\ref{l5.1}), we obtain

\begin{equation}\label{l5.2}
 f(p^q\lambda y)=f((p^{k_1+l_2}c_1d_2-p^{k_2+l_1}c_2d_1)y)
f((p^{l_2}d_2-p^{l_1}d_1)y) f((p^{k_2}c_2-p^{k_1}c_1)y)$$$$
g((p^{k_1+l_2}c_1d_2-p^{k_2+l_1}c_2d_1)y)
g(p^{k_1}c_1(p^{l_2}d_2-p^{l_1}d_1)y)
g(p^{l_1}d_1(p^{k_2}c_2-p^{k_1}c_1)y) $$$$
h((p^{k_1+l_2}c_1d_2-p^{k_2+l_1}c_2d_1)y)
h(p^{k_2}c_2(p^{l_2}d_2-p^{l_1}d_1)y)
h(p^{l_2}d_2(p^{k_2}c_2-p^{k_1}c_1)y)
\end{equation}

Putting  $u=(p^{k_2}c_2-p^{l_2}d_2)y$, $v=(p^{l_2}d_2-I)y$,
$w=(I-p^{k_2}c_2)y$, $y\in Y$, in (\ref{l5.1}), we obtain

\begin{equation}\label{l5.3}
g(p^q\lambda y)=f((p^{l_2}d_2-p^{k_2}c_2)y) f((p^{l_2}d_2-I)y)
f((p^{k_2}c_2-I)y) $$$$ g((p^{l_2}d_2-p^{k_2}c_2)y)
g(p^{k_1}c_1(p^{l_2}d_2-I)y) g(p^{l_1}d_1(p^{k_2}c_2-I)y)
$$$$ h((p^{l_2}d_2-p^{k_2}c_2)y)
h(p^{k_2}c_2(p^{l_2}d_2-I)y) h(p^{l_2}d_2(p^{k_2}c_2-I)y)
\end{equation}

Putting  $u=(p^{l_1}d_1-p^{k_1}c_1)y$, $v=(I-p^{l_1}d_1)y$,
$w=(p^{k_1}c_1-I)y$, $y\in Y$, in (\ref{l5.1}), we obtain
\begin{equation}\label{l5.4} h(p^q\lambda
y)=f((p^{l_1}d_1-p^{k_1}c_1)y) f((p^{l_1}d_1-I)y) f((p^{k_1}c_1-I)y)
$$$$
g((p^{l_1}d_1-p^{k_1}c_1)y) g(p^{k_1}c_1(p^{l_1}d_1-I)y)
g(p^{l_1}d_1(p^{k_1}c_1-I)y)
$$$$
h((p^{l_1}d_1-p^{k_1}c_1)y) h(p^{k_2}c_2(p^{l_1}d_1-I)y)
h(p^{l_2}d_2(p^{k_1}c_1-I)y)
\end{equation}

We check that in fact $q=k$, i.e. the case of $q>k$ is impossible.

Actually, suppose that $q>k$.

Put $y\in Y_{(p^{q})}$ in (\ref{l5.3}). Then $p^q\lambda y=0$, and
the left side of (\ref{l5.3}) is equal to 1. Hence
$g(p^{k}c_1(p^{l_2}d_2-I)y)=1$. Since $p^{l_2}d_2-I \in Aut(Y)$, it
follows from this that $g(y)=1$ for $y\in Y_{(p^{q-k})}$. Let $y\in
Y_{(p^{2q-k})}$ in (\ref{l5.3}). Reasoning similarly we obtain that
$g(y)=1$ for $y\in Y_{(p^{2q-2k})}$. Analogously we obtain that
$g(y)=1$ for $y\in Y_{(p^{N(q-k)})}$ for any natural $N$. Since
$Y=\bigcup_N Y_{(N)}$ and $Y_{(N)}\subset Y_{(N+1)}$, we obtain that
$g(y)=1$ for $y\in Y$. Then (\ref{l5.3}) implies that $f(y)=h(y)=1$
for $y\in Y$. Thus, $\mu_1=\mu_2=\mu_3=E_0$. We obtain the
contradiction to the conditions of the lemma.

Thus, $q=k$.

Note that it follows from equation (\ref{l5.3}) that

\begin{equation}\label{l5.5}
    f(y)=1, y\in Y_{(p^k)}.
\end{equation}
Really, if  $y\in Y_{(p^k)}$, then $g(p^k\lambda y)=1$. Since
$p^{l_2}d_2-I \in Aut(Y)$ and a subgroup $Y_{(p^k)}$ is
characteristic, it follows from(\ref{l5.3}) that equality
(\ref{l5.5}) holds true.

It follows from Lemma 1 that Lemma 7 in cases 1,2,5,8 will be proved
if we show that all solutions of equation (\ref{l5.1}) have form
(\ref{i2}) (the functions are not necessarily equal). In cases 3-4,
6-7, 9-21 of Lemma 7 we construct positive definite solutions
$f(y)$, $g(y)$, $h(y)$ of equation (\ref{l5.1}) such that at least
one of the functions $f(y)$, $g(y)$, $h(y)$ cannot be represented in
the form (\ref{i2}). By the Bochner theorem there exist
distributions $\mu_j$ such that $f(y)=\hat\mu_1(y)$,
$g(y)=\hat\mu_2(y)$, $h(y)=\hat\mu_3(y)$. Then Lemma 1 implies Lemma
7 in cases 3-4, 6-7, 9-21.

We consider all possible cases.

\textbf{Case 1.} $k_1=k_2=l_1=l_2=k$. Since $q=k$, equation
(\ref{l5.3}) takes the form


\begin{equation}\label{l5.1.2}
g(p^k\lambda y)=f(p^{k}(d_2-c_2)y) f((p^{k}d_2-I)y) f((p^{k}c_2-I)y)
$$$$ g(p^{k}(d_2-c_2)y) g(p^{k}c_1(p^{k}d_2-I)y)
g(p^{k}d_1(p^{k}c_2-I)y) $$$$ h(p^{k}(d_2-c_2)y)
h(p^{k}c_2(p^{k}d_2-I)y) h(p^{k}d_2(p^{k}c_2-I)y).
\end{equation}

We get from (\ref{l5.1.2}) that

\begin{equation}\label{l5.1.3}
    g(p^k\lambda y) \leq g(p^{k}c_1(p^{k}d_2-I)y), y \in Y.
\end{equation}
It follows from this that
\begin{equation}\label{l5.1.4}
    g(y) \leq g(\lambda^{-1} c_1(p^{k}d_2-I)y), y \in Y.
\end{equation}
Since each $y\in Y_{(p^{n})}$ for some $n$, $\lambda^{-1}
c_1(p^{k}d_2-I)\in Aut(Y)$, and $Y_{(p^{n})}$ is a characteristic
subgroup, there exist such number $n'$ for every $y$ that
$(\lambda^{-1} c_1(p^{k}d_2-I))^{n'}y=y$. Then it follows from
(\ref{l5.1.4}) that

$$
    g(y) \leq g(\lambda^{-1} c_1(p^{k}d_2-I)y)\leq ... \leq g((\lambda^{-1}
    c_1(p^{k}d_2-I))^{n'}y)=g(y), y \in Y.
$$
So,
$$
    g(y)= g(\lambda^{-1} c_1(p^{k}d_2-I)y)= ... = g((\lambda^{-1}
    c_1(p^{k}d_2-I))^{n'}y)=g(y), y \in Y.
$$
Thus, inequality (\ref{l5.1.3}) becomes an equality

\begin{equation}\label{l5.1.5}
    g(p^k\lambda y) = g(p^{k}c_1(p^{k}d_2-I)y), y \in Y.
\end{equation}
Suppose that $g(p^k\lambda y_0)\neq 0$ for some $y_0\in Y$ such that
$p^k\lambda y_0\neq 0$. Then it follows from (\ref{l5.1.2}) that

\begin{equation}\label{l5.1.6}
1=f(p^{k}(d_2-c_2)y_0) f((p^{k}d_2-I)y_0) f((p^{k}c_2-I)y_0) $$$$
g(p^{k}(d_2-c_2)y_0) g(p^{k}d_1(p^{k}c_2-I)y_0)
$$$$ h(p^{k}(d_2-c_2)y_0)
h(p^{k}c_2(p^{k}d_2-I)y_0) h(p^{k}d_2(p^{k}c_2-I)y_0).
\end{equation}
Since $p^k\lambda y_0\neq 0$, we have $p^{k}d_1(p^{k}c_2-I)y_0\neq
0$ and $p^{k}c_2(p^{k}d_2-I)y_0\neq 0$. So it follows from
(\ref{l5.1.6}) that $g(p^{k}d_1(p^{k}c_2-I)y_0)=1$ and
$h(p^{k}c_2(p^{k}d_2-I)y_0)=1$. By Lemma 5 the set of elements where
a characteristic function is equal to 1 is a subgroup, and each
subgroup of $Y$ contains $Y_{(p)}$, we get that $g(y)=h(y)=1$ for
$y\in Y_{(p)}$. Taking into account (\ref{l5.5}), we get
contradiction to the conditions of the lemma. So, $g(y)=0$ for $y\in
Y\setminus \{0\}$. Then it follows from (\ref{l5.4}) that $h(y)=0$
for $y\in Y\setminus \{0\}$. We get from (\ref{l5.2}) that $f(y)=0$
for $y\in Y\setminus Y_{(p^k)}$.

Thus, we obtain that
\begin{equation}\label{l5.1.7}
    f(y)= \left\{%
\begin{array}{ll}
    1, & \hbox{$y\in Y_{(p^{k})}$;} \\
    0, & \hbox{$y\not\in Y_{(p^{k})}$,} \\
\end{array}%
\right.
g(y)=h(y)= \left\{%
\begin{array}{ll}
   1, & \hbox{$y=0$;} \\
   0, & \hbox{$y\neq 0$.} \\
\end{array}%
\right.
\end{equation}

It follows from (\ref{i2}) and (\ref{l5.1.7}) that $\mu_1=m_{p^k
c_p}$, $\mu_2=\mu_3=m_{c_p}$.

\textbf{Cases 2, 5, 8.} Reasoning as in case 1 we study  the cases
2, 5, 8.

\textbf{Case 3.} $k_1=l_1=k$, $k_2=l_2=l$. Put
\begin{equation}\label{l5.3.9}
    f(y)= \left\{%
\begin{array}{ll}
    1, & \hbox{$y\in Y_{(p^{k+1})}$;} \\
    0, & \hbox{$y\not\in Y_{(p^{k+1})}$,} \\
\end{array}%
\right.
\end{equation}
\begin{equation}\label{l5.3.7}
    g(y)= \left\{%
\begin{array}{ll}
    1, & \hbox{$y\in Y_{(p)}$;} \\
    0, & \hbox{$y\not\in Y_{(p)}$,} \\
\end{array}%
\right.
\end{equation}
\begin{equation}\label{l5.3.8}
    h(y)= \left\{%
\begin{array}{ll}
    h_0(y), & \hbox{$y\in Y_{(p)}$;} \\
    0, & \hbox{$y\not\in Y_{(p)}$,} \\
\end{array}%
\right.
\end{equation}
where $h_0(y)$ is an arbitrary characteristic function on $Y_{(p)}$.
By Lemma 7 the function $h(y)$ is a positive definite function on
$Y$.

We check that the functions (\ref{l5.3.9}), (\ref{l5.3.7}) and
(\ref{l5.3.8}) satisfy equation (\ref{l5.1}).

Let $u\in Y_{(p)}$ in (\ref{l5.1}). Taking into account Lemma 5,
(\ref{l5.3.9}) and (\ref{l5.3.7}), we get

\begin{equation}\label{l5.3.10}
f(v+w) g(p^{k}c_1v+p^{k}d_1w)h(u+p^{l}c_2v+p^{l}d_2w)=$$
$$=f(v)f(w)
g(p^{k}c_1v)g(p^{k}d_1w)h(u)h(p^{l}c_2v)h(p^{l}d_2w), \quad v, w \in
Y.
\end{equation}

Let $v, w\in Y_{(p^{k+1})}$. It follows from Lemma 5 and
(\ref{l5.3.9}) that $f(v+w)=f(v)=f(w)=1$. Since $p^{k}c_1v,
p^{k}d_1w \in Y_{(p)}$, it follows from Lemma 5 and (\ref{l5.3.7})
that $g(p^{k}c_1v+p^{k}d_1w)=g(p^{k}c_1v)=g(p^{k}d_1w)=1$. Since
$p^{l}c_2 v=p^{l}d_2w=0$, we get $h(u+p^{l}c_2v+p^{l}d_2w)=h(u)$,
$h(p^{l}c_2v)=h(p^{l}d_2w)=1$. Thus, equation (\ref{l5.3.10}) holds
true.

Let $v\in Y_{(p^{k+1})}$, $w\not\in Y_{(p^{k+1})}$. Then $p^{k}c_1v
\in Y_{(p)}$ and $p^{k}d_1w \not\in Y_{(p)}$. It follows from Lemma
5 and (\ref{l5.3.7}) that $g(p^{k}c_1v+p^{k}d_1w)=g(p^{k}d_1w)=0$,
and equation (\ref{l5.3.10}) holds true.

Let $v\not\in Y_{(p^{k+1})}$. Then $p^{k}c_1v\not\in Y_{(p)}$ and
the right side of (\ref{l5.3.10}) is equal to 0.

If $w\in Y_{(p^{k+1})}$ then
$g(p^{k}c_1v+p^{k}d_1w)=g(p^{k}c_1v)=0$. So, the left side of
(\ref{l5.3.10}) is also equal to 0. Equation (\ref{l5.3.10}) holds
true.

Let $w\not\in Y_{(p^{k+1})}$. Suppose that the left side of
(\ref{l5.4.2}) does not equal to 0. It follows from (\ref{l5.3.9})
and (\ref{l5.3.7}) that
\begin{equation}\label{l5.3.11}
    v+w \in Y_{(p^{k+1})},
\end{equation}
\begin{equation}\label{l5.3.11.1}
p^{k}c_1v+p^{k}d_1w \in Y_{(p)}.
\end{equation}
It follows from (\ref{l5.3.11}) that $p^k(v+w) \in Y_{(p)}$. Taking
into account (\ref{l5.3.11.1}) this implies that $p^{k}(d_1-c_1)w
\in Y_{(p)}$. Since $q=k$, it is easy to see from (\ref{l5.6}) and
(\ref{l5.7}) that $d_1-c_1\in Aut(Y)$. Therefore we obtain that
$w\in Y_{(p^{k+1})}$. So we obtain the contradiction, and the left
side of equation (\ref{l5.3.10}) is also equal to 0. Thus, we proved
that if $u\in Y_{(p)}$ then the functions (\ref{l5.3.9}),
(\ref{l5.3.7}) and (\ref{l5.3.8}) satisfy equation (\ref{l5.1}).

Let $u\not\in Y_{(p)}$ in (\ref{l5.1}). It follows from
(\ref{l5.3.9}) that $g(u)=0$. The right side of (\ref{l5.1}) is
equal to zero. Suppose that the left side of (\ref{l5.1}) does not
equal to zero. It follows from (\ref{l5.3.9}), (\ref{l5.3.7}) and
(\ref{l5.3.8}) that
\begin{equation}\label{l5.3.12}
    u+v+w  \in Y_{(p^{k+1})},
\end{equation}
\begin{equation}\label{l5.3.12.1}
    u+p^{k}c_1v+p^{k}d_1w \in Y_{(p)},
\end{equation}
\begin{equation}\label{l5.3.12.2}
    u+p^{l}c_2v+p^{l}d_2w \in Y_{(p)}.
\end{equation}
It follows from (\ref{l5.3.12}) that $p^{k}u+p^{k}v+p^{k}w  \in
Y_{(p)}$. Taking into account (\ref{l5.3.12.1}) this implies that
\begin{equation}\label{l5.3.13}
    (p^{k}c_1-I)u+p^{k}(c_1-d_1)w  \in Y_{(p)}.
\end{equation}
It follows from (\ref{l5.3.12}) and the condition $l>k$ that
$p^{l}u+p^{l}v+p^{l}w=0$. Taking into account (\ref{l5.3.12.2}) this
implies that
\begin{equation}\label{l5.3.14}
    (p^{l}c_2-I)u+p^{l}(c_2-d_2)w  \in Y_{(p)}.
\end{equation}
Note that $p^{k}c_1-I, p^{l}c_2-I \in Aut(Y)$. It follows from
(\ref{l5.6}), (\ref{l5.7}), (\ref{l5.3.13}) and (\ref{l5.3.14}) that
$p^k\lambda w\in Y_{(p)}$. Since $\lambda \in Aut(Y)$, we get that
$p^k w\in Y_{(p)}$. Since $c_1-d_1 \in Aut(Y)$, we get that
$p^{k}(c_1-d_1)w  \in Y_{(p)}$, and (\ref{l5.3.13}) implies that
$(p^{k}c_1-I)u\in Y_{(p)}$. Since $p^{k}c_1-I \in Aut(Y)$, we get
that $u\in Y_{(p)}$. So, the obtained contradiction shows that the
left side of equation (\ref{l5.1}) is also equal to 0.

\textbf{Case 4.} $k_1=k_2=k$, $l_1=l_2=l$. Put
\begin{equation}\label{l5.4.1.4}
    f(y)=\left\{%
\begin{array}{ll}
    f_0(y), & \hbox{$y\in Y_{(p^l)}$;} \\
    0, & \hbox{$y\not\in Y_{(p^l)}$,} \\
\end{array}%
\right.
\end{equation}
where $f_0(y)$ is an arbitrary characteristic function on
$Y_{(p^l)}$ such that $f_0(y)=1$ for $y\in Y_{(p^k)}$. By Lemma 7
the function $f(y)$ is a positive definite function on $Y$.

Also put
\begin{equation}\label{l5.4.1.1}
    g(y)=h(y)= \left\{%
\begin{array}{ll}
    1, & \hbox{$y=0$;} \\
    0, & \hbox{$y\neq 0$.} \\
\end{array}%
\right.
\end{equation}

We check that the functions (\ref{l5.4.1.4}) and (\ref{l5.4.1.1})
satisfy equation (\ref{l5.1}).

Let $u=0$ in (\ref{l5.1}). We get

\begin{equation}\label{l5.4.2}
f(v+w) g(p^{k}c_1v+p^{l}d_1w)h(p^{k}c_2v+p^{l}d_2w)=$$
$$=f(v)f(w)
g(p^{k}c_1v)g(p^{l}d_1w)h(p^{k}c_2v)h(p^{l}d_2w), \quad v, w \in Y.
\end{equation}

If $v\in Y_{(p^k)}$ then $p^{k}c_j v=0$ and by Lemma 5 $f(v)=1$,
$f(v+w)=f(w)$. Equation (\ref{l5.4.2}) holds true.

If $v\not\in Y_{(p^k)}$ then $p^{k}c_j v\neq 0$ and
$g(p^{k}c_1v)=0$. Thus, the right side of (\ref{l5.4.2}) is equal to
0.

If $w\in Y_{(p^l)}$ then $p^{l}d_1w=0$ and
$g(p^{k}c_1v+p^{l}d_1w)=g(p^{k}c_1v)=0$. So, the left side is also
equal to 0, and equation (\ref{l5.4.2}) holds true.

Let $w\not\in Y_{(p^l)}$. Suppose that the left side of
(\ref{l5.4.2}) does not equal to 0. It follows from (\ref{l5.4.1.4})
and (\ref{l5.4.1.1}) that
\begin{equation}\label{l5.4.3.1}
    v+w \in Y_{(p^l)},
\end{equation}
\begin{equation}\label{l5.4.3}
    p^{k}c_1v+p^{l}d_1w=0.
\end{equation}
It follows from (\ref{l5.4.3.1}) that $p^l(v+w)=0$. Taking into
account (\ref{l5.4.3}) this implies that $p^k(p^{l-k}d_1-c_1)v=0$.
Since $k<l$, $p^{l-k}d_1-c_1\in Aut(Y)$. Thus, we get that $v\in
Y_{(p^k)}$. So, the obtained contradiction implies that the left
side of equation (\ref{l5.4.2}) is also equal to 0. Thus, we proved
that if $u=0$ then the functions (\ref{l5.4.1.4}) and
(\ref{l5.4.1.1}) satisfy equation (\ref{l5.4.2}).

Let $u\neq 0$ in  (\ref{l5.1}). Then $g(u)=0$ and the right side of
(\ref{l5.1}) is equal to 0. Suppose that the left side of
(\ref{l5.1}) does not equal to 0. Then
\begin{equation}\label{l5.4.4}
    u+v+w \in \mathbb{Z}(p^l),
\end{equation}
\begin{equation}\label{l5.4.4.1}
    u+p^{k}c_1v+p^{l}d_1w=0,
\end{equation}
\begin{equation}\label{l5.4.4.2}
    u+p^{k}c_2v+p^{l}d_2w=0.
\end{equation}
It follows from (\ref{l5.4.4}) that $p^l(u+v+w)=0$. Taking into
account (\ref{l5.4.4.1}) and (\ref{l5.4.4.2}) this imlies that
\begin{equation}\label{l5.4.5}
    (p^l d_1-I)u+(p^ld_1-p^kc_1)v=0, $$$$
    (p^l d_2-I)u+(p^ld_2-p^kc_1)v=0.
\end{equation}
It follows from (\ref{l5.4.5}) that $p^k\lambda v=0$. Hence $v\in
Y{(p^k)}$. Then we obtain from (\ref{l5.4.4}) and (\ref{l5.4.4.1})
that
\begin{equation}\label{l5.4.6}
    p^l(u+w)=0,$$$$
    u+p^{l}d_1w=0.
\end{equation}
It follows from (\ref{l5.4.6}) that $(p^ld_1-I)u=0$. Since
$p^ld_1-I\in Aut(Y)$, we get that $u=0$. So, the obtained
contradiction implies that the left side of equation (\ref{l5.1}) is
also equal to 0. So, the obtained contradiction shows that the left
side of equation (\ref{l5.1}) is also equal to 0.

\textbf{Case 6.} The study of this case is the same as the study of
case 3.

\textbf{Case 7.} The study of this case is the same as the study of
case 4.

\textbf{Case 9.} $k_1=k$, $k_2=l_1=l_2=l$. We check that the
functions (\ref{l5.4.1.4}) and (\ref{l5.4.1.1}) satisfy equation
(\ref{l5.1}).

Let $u=0$ in (\ref{l5.1}). Then equation (\ref{l5.1}) takes the form

\begin{equation}\label{l5.9.16}
f(v+w) g(p^{k}c_1v+p^{l}d_1w)h(p^{l}c_2v+p^{l}d_2w)=$$
$$=f(v)f(w)
g(p^{k}c_1v)g(p^{l}d_1w)h(p^{l}c_2v)h(p^{l}d_2w), \quad v,w\in Y.
\end{equation}

Suppose that either $v \not\in Y{(p^{k})}$, or $w\not\in Y{(p^l)}$.
Then $f(w) g(p^{k}c_1v)=0$, and the right side of (\ref{l5.9.16}) is
equal to 0. Suppose that the left side of (\ref{l5.9.16}) does not
equal to 0. Then
\begin{equation}\label{l5.9.17}
    v+w \in Y_{(p^l)},
\end{equation}
\begin{equation}\label{l5.9.17.1}
    p^{k}c_1v+p^{l}d_1w =0.
\end{equation}
It follows from (\ref{l5.9.17}) that $p^l(v+w)=0$. Taking into
account (\ref{l5.9.17.1}) this implies that
\begin{equation}\label{l5.9.18}
p^{k}(p^{l-k}d_1-c_1)v=0,
$$$$
p^{l}(c_1-p^{l-k}d_1)w=0.
\end{equation}
Since $p^{l-k}d_1-c_1, c_1-p^{l-k}d_1 \in Aut(Y)$, we have $v\in
Y_{(p^{k})}$, $w\in Y_{(p^{l})}$. So, the obtained contradiction
implies that the left side of equation (\ref{l5.9.16}) is also equal
to 0.

Let $u\neq 0$ in (\ref{l5.1}). Then $h(u)=0$, and the right side of
(\ref{l5.1}) is equal to 0


Suppose that the left side of (\ref{l5.1}) does not equal to zero.
Then
\begin{equation}\label{l5.9.20}
    u+v+w\in Y_{(p^l)},$$$$
\end{equation}
\begin{equation}\label{l5.9.20.1}
    u+p^{k}c_1v+p^{l}d_1w=0,$$$$
\end{equation}
\begin{equation}\label{l5.9.20.2}
    u+p^{l}c_2v+p^{l}d_2w=0.
\end{equation}
It follows from (\ref{l5.9.20}) that $u+v+w\in Y_{(p^l)}$.
Taking (\ref{l5.9.20.1}) and (\ref{l5.9.20.2}) into account this
implies that
\begin{equation}\label{l5.9.22}
    (p^{l}d_1-I)u+(p^{l}d_1-p^{k}c_1)v=0,$$$$
    (p^{l}d_2-I)u+(p^{l}d_2-p^{l}c_2))v=0.
\end{equation}
It follows from this that $p^k \lambda v=0$.
Hence $v\in Y_{(p^{k})}$. Since $k<l$ and $p^{l}d_1-I \in Aut(Y)$,
it follows from (\ref{l5.9.22}) that $u=0$. So, the obtained
contradiction implies that the left side of equation (\ref{l5.1}) is
also equal to 0.

\textbf{Case 10.} $k_1=k$, $k_2=l_1=l$, $l_2=m$. We check that the
functions (\ref{l5.4.1.4}) and (\ref{l5.4.1.1}) satisfy equation
(\ref{l5.1}).

Let $u=0$ in (\ref{l5.1}). Then equation (\ref{l5.1}) takes the form

\begin{equation}\label{l5.10.2}
f(v+w) g(p^{k}c_1v+p^{l}d_1w)h(p^{l}c_2v+p^{m}d_2w)=$$
$$=f(v)f(w)
g(p^{k}c_1v)g(p^{l}d_1w)h(p^{l}c_2v)h(p^{m}d_2w), v,w\in Y.
\end{equation}

If $v\in Y_{(p^k)}$ then $p^{k}c_1 v=p^{l}c_2v=0$ and by Lemma 5
$f(v)=1$, $f(v+w)=f(w)$. Equation (\ref{l5.4.2}) holds true.

Let $v\not\in Y_{(p^k)}$. Then $g(p^{k}c_1v)=0$, and a right side of
equation (\ref{l5.10.2}) is equal to 0. Suppose that the left side
of (\ref{l5.10.2}) does not equal to 0. Then
\begin{equation}\label{l5.10.3}
    v+w\in Y_{(p^l)},
\end{equation}
\begin{equation}\label{l5.10.3.1}
    p^{k}c_1v+p^{l}d_1w=0.
\end{equation}
It follows from (\ref{l5.10.3}) that $p^l(v+w)=0$. Taking into
account (\ref{l5.10.3.1}) this implies that
$p^k(p^{l-k}d_1-c_1)v=0$. Since $p^{l-k}d_1-c_1 \in Aut(Y)$, we get
$v\in Y_{(p^k)}$. We obtain the contradiction.

Let $u\neq 0$ in (\ref{l5.1}). Then $g(u)=0$, and the right side of
equation (\ref{l5.1}) is equal to 0. Suppose that the left side of
(\ref{l5.1}) does not equal to zero. Then
\begin{equation}\label{l5.10.5}
    u+v+w\in Y_{(p^l)},
\end{equation}
\begin{equation}\label{l5.10.5.1}
    u+p^{k}c_1v+p^{l}d_1w=0,
\end{equation}
\begin{equation}\label{l5.10.5.2}
    u+p^{l}c_2v+p^{m}d_2w=0.
\end{equation}
It follows from (\ref{l5.10.5}) that $p^l(u+v+w)=0$.
Taking into account (\ref{l5.10.5.1}) and (\ref{l5.10.5.2}) this
implies that
\begin{equation}\label{l5.10.7}
    (p^ld_1-I)u+p^k(p^{l-k}d_1-c_1)v=0,$$$$
    (p^md_2-I)u+p^l(p^{m-l}d_2-c_2)v=0.
\end{equation}
It follows from (\ref{l5.10.7}) that $p^k\lambda v=0$. Hence $v\in
Y_{(p^k)}$. Then it follows from (\ref{l5.10.5.1}) and
(\ref{l5.10.5.2}) that
\begin{equation}\label{l5.10.9}
    u+p^{l}d_1w=0,$$$$
    u+p^{m}d_2w=0.
\end{equation}
It follows from this that $p^l ( d_1-p^{m-l}d_2 )w=0$. Since
$d_1-p^{m-l}d_2 \in Aut(Y)$, we get $w\in Y_{(p^l)}$. Then it
follows from the first equality of (\ref{l5.10.9}) that $u=0$. We
obtain the contradiction. So,  the left side of (\ref{l5.1}) is
equal to zero too.

\textbf{Case 11.} The study of this case is the same as the study of
case 9.

\textbf{Cases 12-21.} The study of these cases is the same as the
study of case 10.

$\blacksquare$

\bigskip

\textbf{Proof of Theorem 1.} Let $Y$ be a character group of $X$. By
Lemma 1 the characteristic functions $\hat\mu_j(y)$ satisfy equation
(\ref{1}). Put $\nu_j=\mu_j*\bar\mu_j$. We have
$\hat\nu_j(y)=|\hat\mu_j(y)|^2 \ge 0, \ j = 1, 2, 3$. It is obvious
that the characteristic functions $\hat\nu_j(y)$ satisfy equation
(\ref{1}) as well.  Hence, when we prove Theorem 1, we may assume
without loss of generality that $\hat\mu_j(y) \ge 0, \ j = 1, 2, 3$,
because $\mu_j$ and $\nu_j$ are either degenerated distributions or
idempotent distributions simultaneously.

Taking into account the form of automorphisms (\ref{t1}), Lemma 2
implies that all $\sigma(\mu_j) \subset B$, where $B= p^n\Delta_p$,
$n\geq 0$, in $X$. Suppose that $\mu_j$ are not degenerate
distributions. We choose the maximum number $N$ such that
$\sigma(\mu_j) \subset B$, where $B= p^N\Delta_p$ in $X$.

Thus, the problem of description of distributions $\mu_1, \mu_2,
\mu_3$ on the group $X = \Omega_p$ is reduces to the problem of
description of distributions $\mu_1, \mu_2, \mu_3$ on the group $X =
\Delta_p$, where the coefficients of the linear forms are
homomorphisms of the group $X = \Delta_p$. Also at least one support
$\sigma(\mu_j)$ is not contained in a proper subgroup of $X$.

Let $X = \Delta_p$. We have $Y \approx \mathbb{Z}(p^{\infty})$.

It follows from Lemma 3 that if all coefficients of the linear forms
are automorphisms of $X = \Delta_p$, then all distributions are
idempotent. It is easy to verify that if $\Lambda$ is not an
automorphism of the group $Y^3$, i.e. $q>0$, then all distributions
$\mu_j$ are degenerated.

So, we will consider the case when at least one of the automorphisms
is not an automorphism of the group $X = \Delta_p$.

Denote $f(y)=\hat\nu_1(y)$, $g(y)=\hat\nu_2(y)$,
$h(y)=\hat\nu_3(y)$. Note that $f(y)\geq 0$, $g(y)\geq 0$, $h(y)\geq
0$. Note that at least one support $\sigma(\mu_j)$ is not contained
in a proper subgroup of $X$ if and only if

\begin{equation}\label{supports}
    \{y\in Y: f(y)=g(y)=h(y)=1\}=\{0\}.
\end{equation}

Lemma 1 implies that the functions $f(y), g(y), h(y)$ satisfy
equation (\ref{1}).

\textbf{I. } Consider the case when $\Lambda=\Lambda_1$. Then
equation (\ref{1}) takes the form
\begin{equation}\label{2}
    f(u + v + w) g(u + \delta_1 v + \varepsilon_1 w)
    h(u + \delta_2 v + \varepsilon_2 w)=$$ $$= f(u)g(u)h(u)
    f(v) g(\delta_1 v) h(\delta_2 v)
    f(w) g(\varepsilon_1 w) h(\varepsilon_2 w), \quad u, v, w \in
    Y.
\end{equation}

Note that
\begin{equation}\label{c1.1.1}
    \det \Lambda=
    \delta_1\varepsilon_2-\delta_2\varepsilon_1-\delta_1+\delta_2+\varepsilon_1-\varepsilon_2.
\end{equation}

Without restricting the generality, we can assume that the
homomorphism $\varepsilon_2$ is not an automorphism of the group
$Y$.



\textbf{I.A.} Suppose that $q=0$, i.e. $\Lambda\in Aut(Y^3)$. In
this case at least one of the homomorphisms $\delta_1, \delta_2,
\varepsilon_1$ is an automorphism.

Suppose that only one of the homomorphism from $\delta_1, \delta_2,
\varepsilon_1$ is an automorphism. Without restricting the
generality, we can assume that $\delta_1 \in Aut(Y)$, $\delta_2,
\varepsilon_1 \not\in Aut(Y)$.

Put

\begin{equation}\label{5}
    f(y)=\left\{%
\begin{array}{ll}
    f_0(y), & \hbox{$y\in Y_{(p)}$;} \\
    0, & \hbox{$y\not\in Y_{(p)}$.} \\
\end{array}%
\right.\quad g(y)=h(y)= \left\{%
\begin{array}{ll}
    1, & \hbox{$y=0$;} \\
    0, & \hbox{$y\neq 0$,} \\
\end{array}%
\right.
\end{equation}
where $f_0(y)$ is an arbitrary characteristic function on $Y_{(p)}$.
By Lemma 7 the function $f(y)$ is a positive definite function on
$Y$ and there exists a distribution $\mu_1$ such that
$\widehat\mu_1(y) = f(y)$. It is clear that $f_0(y)$ can be chosen
in such a way that $\mu_1\not\in I(X)$. Note that $\mu_2=\mu_3=m_X$.

We check that the functions (\ref{5}) satisfy equation (\ref{2}). If
$u=v=0$, then the equation holds true. Let either $u\neq 0$, or
$v\neq 0$. Then $g(u)g(\delta_1 v)=0$ and the right side of equation
(\ref{2})is equal to 0.

Suppose that the left side of (\ref{2}) does not equal to 0. It
follows from (\ref{5}) that
\begin{equation}\label{7}
    u + v + w \in Y_{(p)},
\end{equation}
\begin{equation}\label{7.1}
    u + \delta_1 v + \varepsilon_1 w=0,
\end{equation}
\begin{equation}\label{7.2}
    u + \delta_2 v + \varepsilon_2 w=0.
\end{equation}
Since $\Lambda\in Aut(Y^3)$, it follows from (\ref{7}), (\ref{7.1})
and (\ref{7.2}) that $(u,v,w)\in Y^3_{(p)}$, i.e. $u,v,w\in
Y_{(p)}$. Since $\delta_2, \varepsilon_2\not\in Aut(Y)$, it follows
from (\ref{7.2}) that $u=0$. Since $\varepsilon_1\not\in Aut(Y)$, we
get $\varepsilon_1 w=0$. Taking it into account it follows from
(\ref{7.1}) that $\delta_1 v=0$. Since $\delta_1\in Aut(Y)$, we get
$v=0$. We obtain that $u=v=0$ that contradict to the assumption. So,
the left side of (\ref{2}) is also equal to 0.

Suppose now that at least two of the homomorphisms $\delta_1,
\delta_2, \varepsilon_1$ are automorphisms.

Suppose that $\delta_1, \delta_2 \in Aut(Y)$, $\varepsilon_1 \not\in
Aut(Y)$. Note that in this case $\Lambda$ is an automorphism of the
group $Y^3$ if $\delta_1-\delta_2 \in Aut(Y)$. Then the functions
(\ref{5}) satisfy equation (\ref{2}). The verification is the same
as in previous case.

Suppose that $\delta_1, \varepsilon_1 \in Aut(Y)$, $ \delta_2\not\in
Aut(Y)$. Put

\begin{equation}\label{10}
    f(y)=g(y)=\left\{%
\begin{array}{ll}
    1, & \hbox{$y\in Y_{(p)}$;} \\
    0, & \hbox{$y\not\in Y_{(p)}$.} \\
\end{array}%
\right.\quad h(y)= \left\{%
\begin{array}{ll}
    h_0(y), & \hbox{$y\in Y_{(p)}$;} \\
    0, & \hbox{$y\not\in Y_{(p)}$,} \\
\end{array}%
\right.
\end{equation}
where $h_0(y)$ is an arbitrary characteristic function on $Y_{(p)}$.
Lemma 7 implies that the function $h(y)$ is a positive definite
function on $Y$ and there exists a distribution $\mu_3$ such that
$\widehat\mu_3(y) = h(y)$. It is clear that $h_0(y)$ can be chosen
in such a way that $\mu_3\not\in I(X)$. Note that
$\mu_1=\mu_2=m_{X^{(p)}}$.

We check that functions (\ref{10}) satisfy equation (\ref{2}). If
$u,v,w \in Y_{(p)}$, then $\delta_2 v=\varepsilon_2 w=0$ and
equation (\ref{2}) holds true. Let either $u\not\in Y_{(p)}$, or
$v\not\in Y_{(p)}$, or $w\not\in Y_{(p)}$. Then $f(u)f(v)f(w)=0$ and
the right side of equation (\ref{2})is equal to 0. Suppose that the
left side of (\ref{2}) does not equal to 0. Then
\begin{equation}\label{12}
    \Lambda \left(%
\begin{array}{c}
  u \\
  v \\
  w \\
\end{array}%
\right) \in (Y^3)_{(p)}.
\end{equation}
Since $\Lambda\in Aut(Y^3)$, it follows from (\ref{12}) that
$(u,v,w)\in Y^3_{(p)}$. Hence $u,v,w \in Y_{(p)}$. We obtain the
contradiction. So, the left side of (\ref{2}) is also equal to 0.

Suppose that $\delta_2, \varepsilon_1 \in Aut(Y)$, $ \delta_1\not\in
Aut(Y)$ or $\delta_1, \delta_2, \varepsilon_1 \in Aut(Y)$. Note that
each subgroup of $Y$ is finite and characteristic. If we consider
the restriction of equation (\ref{2}) on any subgroup of $Y$, then
we are in the conditions of Lemma 6. Lemma 6 implies  the statement
of Theorem 1 in this case.

\textbf{I.B.} Suppose that $q>0$, i.e. $\Lambda\not\in Aut(Y^3)$.


Suppose that $\delta_1, \delta_2, \varepsilon_1,\varepsilon_2\not\in
Aut(Y)$. Lemma 9 implies the statement of Theorem 1 in this case.

Suppose now that at least one of the homomorphisms $\delta_1,
\delta_2, \varepsilon_1$ is an automorphism of the group $Y$. Since
$\Lambda\not\in Aut(Y^3)$, it is easy to see that in this case at
least two homomorphisms from $\delta_1, \delta_2, \varepsilon_1$ are
automorphisms.

It is easy to see that $Ker \Lambda \not\subset
Y\times\{0\}\times\{0\}$, $Ker \Lambda \not\subset \{0\}\times
Y\times\{0\}$, $Ker \Lambda \not\subset \{0\}\times\{0\}\times Y$.

At first we suppose that $Ker \Lambda \not\subset \{0\}\times
Y\times Y$. Then there exists an element $(u_0,v_0,w_0) \in Ker
\Lambda$ where $u_0\neq 0$. Put $u=u_0, v=v_0, w=w_0$ in (\ref{2}).
We obtain

\begin{equation}\label{4}
    1= f(u_0)g(u_0)h(u_0)
    f(v_0) g(\delta_1 v_0) h(\delta_2 v_0)
    f(w_0) g(\varepsilon_1 w_0) h(\varepsilon_2 w_0).
\end{equation}

It follows from (\ref{4}) that $f(u_0)=g(u_0)=h(u_0)=1$. Hence
$f(y)=g(y)=h(y)=1$ when $y$ belongs to the subgroup generated by
$u_0$. So we get the contradiction to the condition
(\ref{supports}). Thus, in this case all distributions are
degenerated.

We suppose now that $Ker \Lambda \subset \{0\}\times Y\times Y$.

Then there exists an element $(0,v_0,-v_0) \in Ker \Lambda$, where
$v_0\neq 0$. Put $u=0, v=v_0, w=-v_0$ in (\ref{2}). We obtain

\begin{equation}\label{4.1}
    1= f^2(v_0) g(\delta_1 v_0) h(\delta_2 v_0)
    g(\varepsilon_1 v_0) h(\varepsilon_2 v_0).
\end{equation}

Suppose that $\delta_1, \delta_2, \in Aut(Y)$. It follows from
(\ref{4.1}) that $f(v_0)=g(\delta_1 v_0)=h(\delta_2 v_0)=1$. Note
that any subgroup of $Y$ is characteristic and contains $Y_(p)$.
Taking into account Lemma 5 we get that $f(y)=g(y)=h(y)=1$ on
$Y_(p)$. So we get the contradiction to the condition
(\ref{supports}). Thus, in this case all distributions are
degenerated.

Suppose that $\delta_2, \varepsilon_1 \in Aut(Y)$. It follows from
(\ref{4.1}) that $f(v_0)=g(\varepsilon_1 v_0)=h(\delta_2 v_0)=1$. As
in previous case we get that all distributions are degenerated.

Suppose that $\delta_1, \varepsilon_1\in Aut(Y)$, $\delta_2\not\in
Aut(Y)$. Put $u=(\delta_1\varepsilon_2-\delta_2\varepsilon_1)y$,
$v=(\varepsilon_2-\varepsilon_1)y$, $w=(\delta_2-\delta_1)y$, $y\in
Y$ in  (\ref{2}). We obtain that

\begin{equation}\label{4.2}
f(p^q\lambda
y)=f((\delta_1\varepsilon_2-\delta_2\varepsilon_1)y)
f((\varepsilon_2-\varepsilon_1)y)
f((\delta_2-\delta_1)y)$$$$
g((\delta_1\varepsilon_2-\delta_2\varepsilon_1)y)
g(\delta_1(\varepsilon_2-\varepsilon_1)y)
g(\varepsilon_1(\delta_2-\delta_1)y) $$$$
h((\delta_1\varepsilon_2-\delta_2\varepsilon_1)y)
h(\delta_2(\varepsilon_2-\varepsilon_1)y)
h(\varepsilon_2(\delta_2-\delta_1)y).
\end{equation}

Put $y\in Y_{(p^{q})}$ in (\ref{4.2}). Then $p^q\lambda y=0$, and
the left side of (\ref{4.2}) is equal to 1. Hence
$f((\delta_2-\delta_1)y)=1$. Since $\delta_2-\delta_1 \in Aut(Y)$,
it follows from this that $f(y)=1$ for $y\in Y_{(p^{q})}$. Let $y\in
Y_{(p^{2q})}$ in (\ref{4.2}). Reasoning similarly we obtain that
$f(y)=1$ for $y\in Y_{(p^{2q})}$. Analogously we obtain that
$f(y)=1$ for $y\in Y_{(p^{Nq})}$ for any natural $N$. Since
$Y=\bigcup_N Y_{(N)}$ and $Y_{(N)}\subset Y_{(N+1)}$, we obtain that
$f(y)=1$ for $y\in Y$. Then (\ref{4.2}) implies that $g(y)=h(y)=1$
for $y\in Y$. So we get the contradiction to the condition
(\ref{supports}). Thus, in this case all distributions are
degenerated. We obtain the contradiction to the conditions of the
theorem.

Thus, Theorem 1 is proved for $\Lambda=\Lambda_1$.

\textbf{II.} Consider the case when $\Lambda=\Lambda_2$. Then
equation (\ref{1}) takes the form

\begin{equation}\label{c1}
    f(u + \delta_1 v + w) g(u + v + \varepsilon_1 w)
    h(u + \delta_2 v + \varepsilon_2 w)=$$ $$= f(u)g(u)h(u)
    f(\delta_1 v) g(v) h(\delta_2 v)
    f(w) g(\varepsilon_1 w) h(\varepsilon_2 w), \quad u, v, w \in
    Y.
\end{equation}

Note that
\begin{equation}\label{c1.1}
    \det \Lambda=
    \delta_1\varepsilon_1-\delta_2\varepsilon_1-\delta_1\varepsilon_2+\delta_2+\varepsilon_2-I.
\end{equation}

\textbf{II.A.} Suppose that $q=0$, i.e. $\Lambda\in Aut(Y^3)$.

Let each row and each column of the matrix $\Lambda$ contain at
least two automorphisms of the group $Y$, i.e. there is at least one
automorphism in each couple $(\delta_1,\delta_2)$,
$(\varepsilon_1,\varepsilon_2)$, $(\delta_2,\varepsilon_2)$. Note
that each subgroup of $Y$ is finite and characteristic. If we
consider the restriction of equation (\ref{c1}) to any subgroup of
$Y$, then we are in the conditions of Lemma 6. Lemma 6 implies  the
statement of Theorem 1 in this case.

Let there exist a row or a column of the matrix $\Lambda$ which does
not contain two automorphisms of the group $Y$.

Let $\delta_1,\delta_2 \not\in Aut(Y)$. Note that it follows from
(\ref{c1.1}) that in this case $I-\varepsilon_2 \in Aut(Y)$. Put

\begin{equation}\label{c3}
    \quad f(y)=h(y)= \left\{%
\begin{array}{ll}
    1, & \hbox{$y=0$;} \\
    0, & \hbox{$y\neq 0$,} \\
\end{array}%
\right.
g(y)=\left\{%
\begin{array}{ll}
    g_0(y), & \hbox{$y\in Y_{(p)}$;} \\
    0, & \hbox{$y\not\in Y_{(p)}$,} \\
\end{array}%
\right.
\end{equation}
where $g_0(y)$ is an arbitrary characteristic function on $Y_{(p)}$.
By Lemma 7 the function $g(y)$ is a positive definite function on
$Y$ and there exists a distribution $\mu_2$ such that
$\widehat\mu_2(y) = g(y)$. It is clear that $g_0(y)$ can be chosen
in such a way that $\mu_2\not\in I(X)$. Note that $\mu_1=\mu_3=m_X$.

We check that the functions (\ref{c3}) satisfy equation (\ref{c1}).
If $u=w=0$, then the equation (\ref{c1}) becomes an equality. Let
either $u\neq 0$, or $w\neq 0$. Then $f(u)f(w)=0$ and the right side
of (\ref{c1}) is equal to 0. Suppose that the left side of
(\ref{c1}) does not equal to 0. Then

\begin{equation}\label{c4}
    u + \delta_1 v + w=0,
\end{equation}
\begin{equation}\label{c4.1}
    u +v + \varepsilon_1 w\in Y_{(p)},
\end{equation}
\begin{equation}\label{c4.2}
    u + \delta_2 v + \varepsilon_2 w=0.
\end{equation}
Since $\Lambda\in Aut(Y^3)$, it follows from (\ref{c4}),
(\ref{c4.1}) and (\ref{c4.2}) that $(u,v,w)\in Y^3_{(p)}$. Hence,
$u,v,w\in Y_{(p)}$. Since $\delta_1,\delta_2 \not\in Aut(Y)$, we get
from (\ref{c4}) and (\ref{c4.2}) that
\begin{equation}\label{c4.3}
    u + w=0,
\end{equation}
\begin{equation}\label{c4.4}
    u + \varepsilon_2 w=0.
\end{equation}
It follows from (\ref{c4.3}) and (\ref{c4.4}) that
$(1-\varepsilon_2)w=0$. Since $I-\varepsilon_2 \in Aut(Y)$, we get
that $w=0$ and hence from (\ref{c4.3}) that $u=0$. We obtain the
contradiction to the assumption. So, the left side of (\ref{c1}) is
also equal to 0.

Let $\varepsilon_1,\varepsilon_2 \not\in Aut(Y)$. Note that it
follows from (\ref{c1.1}) that in this case $1-\delta_2 \in Aut(Y)$.
Put
\begin{equation}\label{c5}
    f(y)=\left\{%
\begin{array}{ll}
    f_0(y), & \hbox{$y\in Y_{(p)}$;} \\
    0, & \hbox{$y\not\in Y_{(p)}$.} \\
\end{array}%
\right.\quad g(y)=h(y)= \left\{%
\begin{array}{ll}
    1, & \hbox{$y=0$;} \\
    0, & \hbox{$y\neq 0$,} \\
\end{array}%
\right.
\end{equation}
where $f_0(y)$ is an arbitrary characteristic function on $Y_{(p)}$.
By Lemma 7 the function $f(y)$ is a positive definite function on
$Y$ and there exists a distribution $\mu_1$ such that
$\widehat\mu_1(y) = f(y)$. It is clear that $f_0(y)$ can be chosen
in such a way that $\mu_1\not\in I(X)$. Note that $\mu_2=\mu_3=m_X$.

We check that the functions (\ref{c5}) satisfy equation (\ref{c1}).
If $u=v=0$, then the equation (\ref{c1}) holds true. Let either
$u\neq 0$, or $v\neq 0$. Then the right side of (\ref{c1}) is equal
to 0. Suppose that the left side of (\ref{c1}) does not equal to 0.
Then
\begin{equation}\label{c6}
    u + \delta_1 v + w\in Y_{(p)},
\end{equation}
\begin{equation}\label{c6.1}
    u +v + \varepsilon_1 w=0,
\end{equation}
\begin{equation}\label{c6.2}
    u + \delta_2 v + \varepsilon_2 w=0.
\end{equation}
Since $\Lambda\in Aut(Y^3)$, it follows from (\ref{c6}),
(\ref{c6.1}) and (\ref{c6.2}) that $(u,v,w)\in Y^3_{(p)}$. Hence,
$u,v,w\in Y_{(p)}$. Since $\varepsilon_1,\varepsilon_2 \not\in
Aut(Y)$, we get from (\ref{c6.1}) and (\ref{c6.2}) that
\begin{equation}\label{c6.3}
    u + v=0,
\end{equation}
\begin{equation}\label{c6.4}
    u + \delta_2 v=0.
\end{equation}
It follows from (\ref{c6.3}) and (\ref{c6.4}) that
$(1-\delta_2)v=0$. Since $I-\delta_2 \in Aut(Y)$, we get that $v=0$
and hence from (\ref{c6.3}) that $u=0$. We obtain the contradiction
to the assumption. So, the left side of (\ref{c1}) is also equal to
0.

Let $\delta_2,\varepsilon_2 \not\in Aut(Y)$. Note that it follows
from (\ref{c1.1}) that in this case $1-\delta_1\varepsilon_1 \in
Aut(Y)$. The verification that the functions (\ref{10}) satisfy
equation (\ref{c1}) is the same as in case I.A ($\delta_1,
\varepsilon_1 \in Aut(Y)$, $ \delta_2\not\in Aut(Y)$).

\textbf{II.B.} Suppose that $q>0$, i.e. $\Lambda\not\in Aut(Y^3)$.
Note that it follows from (\ref{c1.1}) that in this case either
$\delta_1\varepsilon_1\in Aut(Y)$, or $\delta_2\in Aut(Y)$, or
$\varepsilon_2\in Aut(Y)$.

Suppose that either $\delta_2\in Aut(Y)$, or $\varepsilon_2\in
Aut(Y)$. It is easy to see that $Ker \Lambda \not\subset
Y\times\{0\}\times\{0\}$, $Ker \Lambda \not\subset \{0\}\times
Y\times\{0\}$, $Ker \Lambda \not\subset \{0\}\times\{0\}\times Y$.
Hence, there exists an element $(u_0,v_0,w_0)\in Ker \Lambda$ such
that at least two coordinates do not equal to 0. Put $u=u_0, v=v_0,
w=w_0$ in (\ref{c1}). We obtain

\begin{equation}\label{c9}
    1= f(u_0)g(u_0)h(u_0)
    f(\delta_1 v_0) g(v_0) h(\delta_2 v_0)
    f(w_0) g(\varepsilon_1 w_0) h(\varepsilon_2 w_0).
\end{equation}
If $u_0\neq 0$, we get that $f(y)=g(y)=h(y)=1$ when $y$ belongs a
subgroup generated by $u_0$. So we get the contradiction to the
condition (\ref{supports}).

Let $u_0=0$. Then it follows from (\ref{c9}) that $1=
    f(\delta_1 v_0) g(v_0) h(\delta_2 v_0)
    f(w_0) g(\varepsilon_1 w_0) h(\varepsilon_2 w_0).$
Since each subgroup of $Y$ contains $Y_{(p)}$ and either
$\delta_2\in Aut(Y)$, or $\varepsilon_2\in Aut(Y)$, Lemma 5 implies
that $f(y)=g(y)=h(y)=1$ on $Y_{(p)}$. So we get the contradiction to
the condition (\ref{supports}). Thus in these cases all
distributions are degenerated.

Suppose that $\delta_1\varepsilon_1\in Aut(Y)$, $\delta_2\not\in
Aut(Y)$, $\varepsilon_2\not\in Aut(Y)$.

Put $u=(\varepsilon_2-\delta_2\varepsilon_1)y$,
$v=(\varepsilon_1-\varepsilon_2)y$, $w=(\delta_2-I)y$, $y\in Y$ in
(\ref{c1}). We obtain that

\begin{equation}\label{c10}
f(p^q\lambda y)=f((\varepsilon_2-\delta_2\varepsilon_1)y)
f(\delta_1(\varepsilon_1-\varepsilon_2)y) f((\delta_2-I)y)$$$$
g((\varepsilon_2-\delta_2\varepsilon_1)y)
g((\varepsilon_1-\varepsilon_2)y) g(\varepsilon_1(\delta_2-I)y) $$$$
h((\varepsilon_2-\delta_2\varepsilon_1)y)
h(\delta_2(\varepsilon_1-\varepsilon_2)y)
h(\varepsilon_2(\delta_2-I)y).
\end{equation}

Put $y\in Y_{(p^{q})}$ in (\ref{c10}). Then $p^q\lambda y=0$, and
the left side of (\ref{c10}) is equal to 1. Hence
$f((\delta_2-I)y)=1$. Since $\delta_2-I \in Aut(Y)$, it follows from
this that $f(y)=1$ for $y\in Y_{(p^{q})}$. Let $y\in Y_{(p^{2q})}$
in (\ref{c10}). Reasoning similarly we obtain that $f(y)=1$ for
$y\in Y_{(p^{2q})}$. Analogously we obtain that $f(y)=1$ for $y\in
Y_{(p^{Nq})}$ for any natural $N$. Since $Y=\bigcup_N Y_{(N)}$, we
obtain that $f(y)=1$ for $y\in Y$. Then (\ref{c10}) implies that
$g(y)=h(y)=1$ for $y\in Y$. So we get the contradiction to the
condition on supports of $\mu_j$. Thus, in this case all
distributions are degenerated. We obtain the contradiction to the
conditions of the theorem.

$\blacksquare$

\bigskip

\textbf{Remark 2.} Statement А of Theorem 1 can not be strengthened
to the statement that all distributions $\mu_1, \mu_2, \mu_3$ are
degenerated. Namely, it is easy to verify that if $\mu_1=
\mu_2=\mu_3=m_{\Delta_p}$ then components of the random vector
$(L_1,L_2,L_3)$ are independent.

Statement В.2 of Theorem 1 in case $q=k$ and or $k_1=l_2=k$, or
$k_2=l_1=k$ also can not be strengthened to the statement that all
distributions $\mu_1, \mu_2, \mu_3$ are degenerate. , які
характеризуються незалежністю компонент випадкового вектору
$(L_1,L_2,L_3)$. Namely, it is easy to verify that if $\mu_1= m_{p^k
\Delta_p}, \mu_2=\mu_3=m_{\Delta_p}$, then components of the random
vector $(L_1,L_2,L_3)$ are independent.

\end{document}